\documentclass [12pt] {article}
\usepackage{amsfonts}
\newcommand{\qed} {\hspace {0.1in} \rule {1.5mm} {3.5mm}}

\newtheorem{conjecture}{Conjecture}
\newtheorem{theorem}{Theorem}
\newtheorem{proposition}{Proposition}

\def\mk{M^k_{\alpha}}
\def\kmk{\Ke\,\mk}
\def\rmk{\Ra\,\mk}
\def\df{\dim_{\cf_k}}
\def\dg{\dim_G}

\def\sq{\subseteq}

\def\dim{{\rm dim}}
\def\supp{{\rm supp}}
\def\Ke{{\rm Ker}}
\def\Ra{{\rm Ran}}

\def\<{\langle}
\def\>{\rangle}
\def\proof{\smallskip\noindent{\bf Proof:} }

\def\bN{{\mathbb N}}
\def\bC{{\mathbb C}}

\def\cf{\mbox{$\cal F$}}
\def\cg{\mbox{$\cal G$}}

\def\to{\rightarrow}

\title{On the analytic zero divisor conjecture of Linnell}
\author{{\sc G\'abor Elek}
\cr Mathematical Institute of
the Hungarian Academy of Sciences\cr P.O. Box 127, H-1364 Budapest, Hungary\cr
elek@renyi.hu}
\date{}
\begin{document}

\maketitle
\vskip 0.2in
\noindent{\bf Abstract.} In this note we prove that in the case of finitely
generated amenable groups the classical zero divisor conjecture implies the
analytic zero divisor conjecture of Linnell.
\newpage
\section{Introduction}
Let $G$ be a discrete group. 
Denote by $\bC G$ the complex group algebra of $G$. 
The following conjecture
is called the zero divisor conjecture.
\begin{conjecture}
Let $G$ be a torsion free group.
If $0\neq\alpha\in\bC G$ and $0\neq\beta\in\bC G$, then $\alpha\beta\neq 0$.
\end{conjecture}
Now let
$$
L^2(G)=\{\alpha:\,G\to\bC\,\mid\,\,\sum_{g\in G}|\alpha(g)|^2<\infty\,\}\,.
$$
Then the Hilbert space $L^2(G)$ is a two-sided $\bC G$-module. 
In \cite{Lin},
Linnell formulated an analytic version of the zero divisor conjecture.
\begin{conjecture}
Let $G$ be a torsion free group.
If $0\neq\alpha\in\bC G$ and $0\neq\beta\in L^2(G)$, then $\alpha\beta\neq 0$.
\end{conjecture}
Since $\bC G\subseteq L^2(G)$, the second conjecture implies the first one.
The goal of this paper is to prove that for finitely generated amenable groups
the two conjectures are actually equivalent.

\begin{theorem}
Let $G$ be a finitely generated amenable group. 
If $0\neq\alpha\in \bC G$,
$0\neq\beta\in L^2(G)$ and $\alpha\beta=0$, then there exists $0\neq\gamma\in
\bC G$ such that $\alpha\gamma=0$.
\end{theorem}
\section{Invariant subspaces and the von Neumann dimension}

A linear subspace $V\sq L^2(G)$ is called
an invariant subspace it is a right $\bC G$-submodule of $L^2(G)$, that is if $\beta\in V$ and $\gamma\in
\bC G$, then $\beta\gamma\in V$. The von Neumann dimension of an invariant
subspace $V$ is a real number defined the following way \cite{Pan},
$$
\dim_GV=\langle P_V 1_e,1_e\rangle\,
$$
where $P_V$ denotes the orthogonal projection onto the closure of $V$ and
$1_e\in L^2(G)$ takes the value $1$ on the unit of $G$ and vanishes everywhere
else. 
We list some properties of the von Neumann dimension \cite{Pan}.
\begin{enumerate}
\item If $V\neq 0$, then $\dim_G\, (V) > 0$\,.
\item $\dim_G\,( L^2(G))=1$\,.
\item If $V\sq W$ are invariant subspaces, then $\dim_G (V)\leq \dim_G (W)$.
\item If $T:\,L^2(G)\to L^2(G)$ is a bounded linear transformation and a right
$\bC G$-module homomorphism, then

$$
\dim_G (\Ke\,T)+\dim_G (\Ra\,T)=1\,.$$\end{enumerate}
For an arbitrary linear subspace $W\sq L^2(G)$ and a finite subset $A\sq G$
one can consider
$$
\dim_A (W)=\frac{\sum_{g\in A} \langle P_W 1_g, 1_g\rangle} {|A|}\,,
$$
where $1_g\in L^2(G)$  takes the value $1$ on $g$ and vanishes everywhere
else.  \cite{Eck}.
Obviously, $\dim_A (W)$ is a non-negative real number, and if
$W\sq Z$, then $\dim_A(W)\leq\dim_A(Z)$. If $V$ is an invariant
subspace, then for any finite subset $A\sq G$, $\dim_A(V)=\dim_G(V)$. Note
that
if $W$ is a finite dimensional subspace consisting of elements supported on
$A$, then $\dim_A\,(W)=\frac{\dim_{\bC}\,W}{|A|}$\,.
\section{Approximation of the von Neumann dimension}
The goal of this section is to prove a variant of a result of Dodziuk and
Mathai \cite{DM}. First let us recall the notion of amenability.
Let $G$ be a finitely generated group with symmetric generator set
$\{g_1,g_2,\dots,g_k\}$. The generators determine a word-metric $d$ on $G$
such a way that the right multiplications by the elements of $G$ are
isometries. 
The group $G$ is called amenable if it has an exhaustion
$$\cf_1\sq\cf_2\sq\dots, \bigcup_{k=1}^\infty\cf_k=G\,$$
by finite subsets such that for any fixed $r\in \bN$,
$$\lim_{k\to\infty}\frac{|\partial_r\cf_k|}{|\cf_k|}=0\,,$$
where $\partial_r\cf_k=\{g\in\cf_k\,|\,d(g,G\backslash\cf_k)\leq r\}\,.$
Before stating our proposition we list some notations and definitions.
Let $\alpha\in\bC G$, then $M_{\alpha}:\,L^2(G)\to L^2(G)$ is
the left multiplication by $\alpha$. Then $\Ke\,M_{\alpha}$ is an invariant
subspace.
The width of $\alpha$, $w(\alpha)$ is defined as $\max_{g\in G,\,\alpha(g)\neq
0} d(e,g)$, where $e$ is the unit of $G$. Let $\cg_k=\cf_k\backslash
\partial_{w(\alpha)}\cf_k$.
Note that if $\beta$ is supported on $\cg_k$, then $\alpha\beta$ will be
supported on $\cf_k$. We consider the following spaces,
$$
V_k=\{\beta\in L^2(G)\,|\,\supp\,(\beta)\sq\cf_k\}
$$
$$
W_k=\{\beta\in L^2(G)\,|\,\supp\,(\beta)\sq\cg_k\}
$$
\begin{proposition}
$\dim_G(\Ke\, M_\alpha)=\lim_{k\to\infty}\dim_{\cf_k}
(V_k\cap \,\Ke\,M_{\alpha})$\,.
\end{proposition}
\proof
Let $\mk:\,W_k\to V_k$ the
restriction of $M_\alpha$ onto $W_k$.
Then $\kmk\sq V_k\cap\Ke\,M_\alpha$, $\rmk\sq\Ra\, M_\alpha$ and
$\dim_{\bC} \kmk + \dim_{\bC} \rmk=|\cg_k|\,$.
Therefore,
\begin{equation}\label{1e}
\df(\kmk)\leq \df(V_k\cap \Ke\,M_\alpha)\leq \dg(\Ke\,M_\alpha)
\end{equation}
\begin{equation}\label{2e}
\df(\rmk)\leq\dg(\Ra\,M_\alpha)
\end{equation}
\begin{equation}\label{3e}
\df(\kmk)+\df(\rmk)=\frac{|\cg_k|}{|\cf_k|}
\end{equation}
By (\ref{3e}),
\begin{equation}\label{4e}
\lim_{k\to\infty}(\df(\kmk)+\df(\rmk))=1=\dg (\Ke\,M_\alpha)+
\dg(\Ra\,M_\alpha)
\end{equation}
Hence, by (\ref{1e}) and (\ref{2e}), $\lim_{k\to\infty}\df\kmk=
\dg(\Ke\,M_\alpha)$, thus the Proposition follows.\qed

\section{The proof of the Theorem}
Let $0\neq\alpha\in\bC G$, $0\neq\beta\in L^2(G)$ and $\alpha\beta=0$.
Then $\Ke\,M_\alpha\neq 0$, hence $\dg (\Ke\,M_\alpha)>0$. By the Proposition,
for sufficiently large $k$, there exists $0\neq\gamma\in V_k\cap
\Ke\,M_\alpha$\,. Then $\gamma\in\bC G$ and $\alpha\gamma=0$.\qed

\end{document}